\documentclass[psamsfont,a4paper,11pt]{amsart}

\usepackage[english]{babel}   

\usepackage{setspace}     
\usepackage{enumerate} 	
\usepackage{afterpage} 	
\usepackage[usenames,dvipsnames]{color}
\usepackage{graphicx}	

\usepackage{amsmath} 	
\usepackage{amssymb}	
\usepackage{mathrsfs} 	
\usepackage{amsfonts} 	
\usepackage{latexsym} 	
\usepackage[latin1]{inputenc} 

\usepackage{amsthm} 	
\usepackage{stackrel}       
\usepackage{amscd} 	

\usepackage{tabularx}	
\usepackage{longtable}	

\usepackage{verbatim} 
\usepackage{blindtext} 

\allowdisplaybreaks

\def\a{\alpha} 
 
\def\b{\beta} 
\def\d{\delta} 
 
\def\g{\gamma} 
 
\def\l{\lambda}

\def\cal{\mathcal}

\def\H{{\cal H}}

\newcommand{\C}{\mathbb{C}}

\newcommand{\R}{\mathbb{R}}

\newcommand{\Z}{\mathbb{Z}} 
\newcommand{\Q}{\mathbb{Q}}

\renewcommand{\to}{\longrightarrow}

\def\mod{\mathop{\rm{mod}}}

\newtheorem{Thm}{Theorem}[]		

\newtheorem*{thm}{Theorem}	

\newtheorem*{prop}{Proposition}

\theoremstyle{definition}

\newtheorem*{defn}{Definition}  

\theoremstyle{remark}
 
\newtheorem*{rmk}{Remark}

\newcommand{\blip}{\vspace{0.2cm}}

\newtheorem{ind}[]{{\rm\it Indice}}

\usepackage{multicol}

\usepackage{tikz}
\usetikzlibrary{arrows}

\usepackage{paralist}
\usepackage{cite}

\title{Congruences for Powers of the Partition Function}

\author[Locus]{Madeline Locus}
\author[Wagner]{Ian Wagner}

\begin{document}

\maketitle

\begin{abstract}
Let $p_{-t}(n)$ denote the number of partitions of $n$ into $t$ colors.  
In analogy with Ramanujan's work on the partition function,  Lin recently proved in \cite{Lin} that $p_{-3}(11n+7)\equiv0\pmod{11}$ for every integer $n$. Such congruences, those of the form
$p_{-t}(\ell n + a) \equiv 0 \pmod {\ell}$, were previously studied by Kiming and Olsson.  If $\ell \geq 5$ is prime and $-t \not \in \{\ell - 1, \ell -3\}$, then such congruences
satisfy $24a \equiv -t \pmod {\ell}$.  Inspired by Lin's example,
we obtain natural infinite families of such congruences.  If $\ell\equiv2\pmod{3}$ (resp. $\ell\equiv3\pmod{4}$ and $\ell\equiv11\pmod{12}$) is prime and $r\in\{4,8,14\}$ (resp. $r\in\{6,10\}$ and $r=26$), then for $t=\ell s-r$, where $s\geq0$, we have that
\begin{equation*}
p_{-t}\left(\ell n+\frac{r(\ell^2-1)}{24}-\ell\Big\lfloor\frac{r(\ell^2-1)}{24\ell}\Big\rfloor\right)\equiv0\pmod{\ell}.
\end{equation*}
Moreover, we exhibit infinite families where such congruences cannot hold.

\end{abstract}

\section{Introduction and Statement of Results}

A $t$-color partition, where $t \geq 1$, of an integer $n$ is a $t$-tuple $\lambda = (\lambda^{(1)}, ..., \lambda^{(t)})$ of partitions such that $|\lambda^{(1)}| + ... + |\lambda^{(t)}| = n$ (see \cite{Andrews}).  The $t$-color partition function, which counts partitions of $n$ into $t$ colors, has the following generating function:
\[\sum_{n=0}^\infty p_{-t}(n)q^n = \prod_{n=1}^\infty \frac{1}{(1-q^n)^t}.\]
When $t=1$, this becomes the usual partition function $p(n)$.  Ramanujan discovered the following congruences for the partition function:
\[p(5n + 4) \equiv 0 \pmod 5,\]
\[p(7n + 5) \equiv 0 \pmod 7,\]
\[p(11n + 6) \equiv 0 \pmod {11}.\]\\
Notice that the modulus of the progression matches the modulus of the congruence in these three examples.
Motivated by this observation, it is natural to classify other congruences where these moduli agree.  Kiming and Olsson proved in \cite{KO} that in order for a Ramanujan-type congruence modulo $\ell$ to exist for the progression $\ell n+a$, it must be true that $24a\equiv-t\pmod{\ell}$.  
Along these lines, there have been other works which study specific cases in which congruences of the partition function exist.  In \cite{Harper}, Harper gave several arithmetic progressions of primes and stated the possible values of $b$ such that the coefficients $a(n)$ of certain color partitions satisfy similar congruences, namely $a(\ell n+b)\equiv0\pmod{\ell}$.  More recently, in \cite{Lin} Lin found and proved the specific congruence $p_{-3}(11n+7)\equiv0\pmod{11}$ using elementary methods.  In view of these recent works, it is natural to ask if there is a general phenomenon of when congruences can and cannot occur.  Is there a more general form of Harper's and Lin's work?  Kiming and Olsson raise the question of whether there are infinitely many primes for which a congruence of this type exists.  We address both of these questions in the following theorem. 

\begin{Thm}\label{main}
Assume the notation above, and let $t=\ell s-r$.  Then the following are true.
\begin{enumerate}
\item If $r \in \{4, 8, 14\}$ and $\ell \equiv 2 \pmod 3$ is prime, then for every integer $n$, we have
\begin{equation*}
p_{-t} \left( \ell n + \frac{r(\ell^{2}-1)}{24} - \ell\Big\lfloor\frac{r(\ell^2-1)}{24\ell}\Big\rfloor \right) \equiv 0 \pmod \ell.
\end{equation*}
\item If $r \in \{4, 8, 14\}$ and $\ell\equiv1\pmod{3}$, then there are infinitely many $n$ for every arithmetic progression mod $\ell$ for which
\begin{equation*}
p_{-t} \left(\ell n + \frac{r(\ell^{2}-1)}{24} -\ell\Big\lfloor\frac{r(\ell^2-1)}{24\ell}\Big\rfloor \right)\not\equiv0\pmod{\ell}
\end{equation*}
\end{enumerate}
\end{Thm}

\begin{Thm} \label{two}
Assume the notation above.  Then the following are true.
\begin{enumerate}
\item If $r \in \{6, 10\}$ and $\ell \equiv 3 \pmod 4$ is prime, then for every integer $n$, we have
\begin{equation*}
p_{-t} \left(\ell n + \frac{r(\ell^{2}-1)}{24} -\ell\Big\lfloor\frac{r(\ell^2-1)}{24\ell}\Big\rfloor \right) \equiv 0 \pmod \ell.
\end{equation*}
\item If $r \in \{6, 10\}$ and $\ell\equiv1\pmod{4}$, then there are infinitely many $n$ for every arithmetic progression mod $\ell$ for which
\begin{equation*}
p_{-t} \left(\ell n + \frac{r(\ell^{2}-1)}{24} -\ell\Big\lfloor\frac{r(\ell^2-1)}{24\ell}\Big\rfloor \right)\not\equiv0\pmod{\ell}
\end{equation*}
\end{enumerate}
\end{Thm}

\begin{Thm} \label{three}
Assume the notation above.  Then the following are true.
\begin{enumerate}
\item If $r=26$ and $\ell \equiv 11 \pmod {12}$ is prime, then for every integer $n$, we have
\begin{equation*}
p_{-t} \left(\ell n + \frac{r(\ell^{2}-1)}{24} -\ell\Big\lfloor\frac{r(\ell^2-1)}{24\ell}\Big\rfloor \right) \equiv 0 \pmod \ell.
\end{equation*}
\item If $r=26$ and $\ell\equiv1, 5, 7\pmod{12}$, then there are infinitely many $n$ for every arithmetic progression mod $\ell$ for which
\begin{equation*}
p_{-t} \left(\ell n + \frac{r(\ell^{2}-1)}{24} -\ell\Big\lfloor\frac{r(\ell^2-1)}{24\ell}\Big\rfloor \right)\not\equiv0\pmod{\ell}.
\end{equation*}
\end{enumerate}
\end{Thm}

\begin{rmk}
Results similar to Theorems 1, 2, and 3 were already known from a work of Boylan \cite{Boylan}.
\end{rmk}

To prove these statements, we use a famous work of Serre on powers of Dedekind's eta-function (see \cite{Serre}), together with the theory of Hecke operators, modular forms with complex multiplication, and some combinatorial observations.

\section{Nuts and Bolts}

In Section 2.1, we give basic definitions from the theory of modular forms and Hecke operators, and we introduce a theorem of Serre which gives the basis for our proofs of Theorems 1, 2, and 3.  In Section 2.2, we define newforms with complex multiplication, which we reveal in Section 3 are closely connected with Ramanujan-type congruences.

\subsection{Basic Definitions}

Recall the \emph{congruence subgroup} $\Gamma_0(N)$ defined as
\begin{eqnarray*}
\Gamma_0(N)&:=&\left\{\begin{pmatrix}a&b\\c&d\end{pmatrix}\in\mathrm{SL}_2(\Z):c\equiv0\pmod{N}\right\}.
\end{eqnarray*}

Suppose that $f(z)$ is a meromorphic function on $\H$, the upper half of the complex plane, that $k\in\Z$, and that $\Gamma$ is a congruence subgroup of level $N$. Then $f(z)$ is called a \emph{meromorphic modular form with integer weight $k$ on $\Gamma$} if the following hold:
\begin{enumerate}[(1)]
\item We have
\begin{equation*}
f\left(\frac{az+b}{cz+d}\right)=(cz+d)^kf(z)
\end{equation*}
for all $z\in\H$ and all $\begin{pmatrix}a&b\\c&d\end{pmatrix}\in\Gamma$.
\item If $\g_0\in\mathrm{SL}_2(\Z)$, then $(f\vert_k\g_0)(z)$ has a Fourier expansion of the form
\begin{equation*}
\left(f\vert_k\g_0\right)(z)=\sum_{n\geq n_{\g_0}}a_{\g_0}(n)q_N^n,
\end{equation*}
where $q_N:=e^{2\pi iz/N}$ and $a_{\g_0}(n_{\g_0})\neq0$.
If $k=0$, then $f(z)$ is known as a \emph{modular function on $\Gamma$}.
\end{enumerate}
Furthermore, if $\chi$ is a Dirichlet character modulo $N$, then we say that a form $f(z)\in M_k(\Gamma_1(N))$ (resp. $S_k(\Gamma_1(N))$\,) has \emph{Nebentypus character $\chi$} if
\begin{equation*}
f\left(\frac{az+b}{cz+d}\right)=\chi(d)(cz+d)^kf(z)
\end{equation*}
for all $z\in\H$ and all $\begin{pmatrix}a&b\\c&d\end{pmatrix}\in\Gamma_0(N)$.

The space of such modular forms (resp. cusp forms) is denoted by $M_k(\Gamma_0(N),\chi)$ (resp. $S_k(\Gamma_0(N),\chi)$\,).  If $f(z)$ is an integer weight meromorphic modular form on a congruence subgroup $\Gamma$, then we call $f(z)$ a \emph{holomorphic modular (resp. cusp) form} if $f(z)$ is holomorphic on $\H$ and is holomorphic (resp. vanishes) at the cusps of $\Gamma$.

We now define the Hecke operators, which are linear transformations that act on spaces of modular forms by transforming the coefficients of their Fourier expansions.  Let $k\in\Z^+$ and let $p$ be prime.  If $f(z)=\sum_{n=1}^{\infty}A(n)q^n$ is in $M_k(\Gamma_0(N),\chi)$ (resp. $S_k(\Gamma_0(N),\chi)$), then the action of the \emph{Hecke operator} $T_{p,k,\chi}$ on $f(z)$ is defined by
\begin{equation*}
f(z)\mid T_{p,k,\chi}:=\sum\limits_{n=0}^{\infty}\left(A(pn)+\chi(p)p^{k-1}A\left(\frac{n}{p}\right)\right)q^n.
\end{equation*}
If $p\nmid n$, then we say that $A(n/p)=0$.  For our purposes, we shorten $T_{p,k,\chi}$ to $T_p$.  It will later be valuable to note that $f(z)\mid T_{p,k,\chi} \equiv 0 \pmod p$ only when $A(pn) \equiv 0 \pmod p$ for every $n$.

We define the following function as a power of a Dedekind eta-product, with $r,\d\in\Z$.  We only consider integer weight modular forms, so we assume that $r$ is even.
\begin{equation*}
F_{r, \d}(z):=\eta^r(\d z)=q^\frac{\d r}{24}\prod\limits_{n=1}^{\infty}(1-q^{\d n})^r = \sum_{n=0}^\infty A_{r}(n)q^n \in S_k(\Gamma_0(N),\chi). \tag{*}
\end{equation*}
When $r\in\{2,4,6,8,10,14,26\}$, the level $N$ of $\eta^r(\d z)$ is the corresponding value of $\d^2$ in Theorem \ref{Serre}, and the Nebentypus $\chi$ depends on $N$.
Several nice properties of this family of functions will become clear when we introduce newforms in the next section.

We now give a result of Serre from \cite{Serre}, which follows from the theory of modular forms with complex multiplication.
\begin{Thm}[Serre] \label{Serre}
Suppose that $r\in\{2,4,6,8,10,14,26\}$.  Then the following are true:
\begin{enumerate}
\item If $\ell\equiv2\pmod{3}$ is prime, $r\in\{4,8,14\}$, and $\d$ is the corresponding $\{6, 3, 12\}$, then $F_{r, \d}(z)\vert T_\ell \equiv 0 \pmod {\ell}$.
\item If $\ell\equiv3\pmod{4}$ is prime, $r\in\{6,10\}$, and $\d$ is the corresponding $\{4, 12\}$, then\\ $F_{r, \d}(z)\vert T_\ell \equiv 0 \pmod {\ell}$.
\item If $\ell \equiv 11 \pmod {12}$ is prime, $r=26$, and $\d = 12$, then $F_{r, \d}(z) \vert T_{\ell} \equiv 0 \pmod {\ell}$.\\
\end{enumerate}
\end{Thm}

\subsection{Modular forms with Complex Multiplication}

For a generic number field $K$ which is a finite field extension of $\Q$, we denote the ring of algebraic integers in $K$ by $\mathcal{O}_K$.  This ring is a finite field extension of $\Z$, and in the cases we will consider will be a principal ideal domain.

\begin{rmk}
In general, $\mathcal{O}_K$ is rarely a unique factorization domain, but for the number fields we consider, its ideals always factor uniquely into products of prime ideals since it is a principal ideal domain.
\end{rmk}

It will be beneficial to consider $(p)=(a+bi)(a-bi)$ as ideals in $\Z[i]$, the ring of integers of $\mathcal{O}_K$ where $K=\Q(i)$, and $(p) = \left(x+y \sqrt{-3}\right)\left(x-y \sqrt{-3}\right)$ as ideals in $\Z\big[\frac{1 + \sqrt{-3}}{2}\big]$, the ring of integers of $\mathcal{O}$ where $K = \Q(\sqrt{-3})$, for those primes which factor (i.e. ``the split primes"). These are the unique factorizations of $(p)$ into prime ideals in their respective fields.  In order to connect this object with congruences for powers of the partition function we will make use of newforms.
A \emph{newform} is a normalized cusp form in $S_k^{\text{new}}\left(\Gamma_0(N)\right)$ that is an eigenform of all the Hecke operators and all the Atkin-Lehner involutions $\vert_kW(Q_p)$, for primes $p\vert N$, and $\vert_kW(N)$.

Newforms come with the following properties.

\begin{thm}
If $f(z)=\sum\limits_{n=1}^\infty A_r(n)q^n\in S_k^{\text{new}}\left(\Gamma_0(N)\right)$ is a newform, then the following are true.
\begin{enumerate}
\item If $p$ is prime with $p^2\vert N$, then $A_r(p)=0$.
\blip
\item If $p\vert N$ is prime with $\mathrm{ord}_p(N)=1$, then $A_r(p)=-\l_pp^{\frac{k}{2}-1}$, where $\l_p$ is a root of unity.
\end{enumerate}
\end{thm}

If $K = \Q(\sqrt{-D})$ is an imaginary quadratic field with discriminant $-D$ and $\mathcal{O}_{K}$ is its ring of algebraic integers, then a Hecke Grossencharacter $\phi$ of weight $k\geq 2$ with modulus $\Lambda$ can be defined in the following way.  Let $\Lambda$ be a nontrivial ideal in $\mathcal{O}_{K}$ and let $I(\Lambda)$ denote the group of fractional ideals prime to $\Lambda$.  A Hecke Grossencharacter $\phi$ with modulus $\Lambda$ is a homomorphism 
\[\phi : I(\Lambda) \rightarrow \C^{\times} \]
such that for each $\alpha \in K^{\times}$ with $\alpha \equiv 1 \pmod {\Lambda}$ we have 
\[ \phi(\alpha \mathcal{O}_{K}) = \alpha^{k-1}. \]
Using this notation the newform $\Psi(z)$ can be defined by
\[\Psi(z) = \sum_{(\alpha)} \phi((\alpha))q^{N((\alpha))} = \sum_{n=1}^\infty a(n) q^n. \]
The sum is over the integral ideals $(\alpha)$ that are prime to $\Lambda$ and $N((\alpha))$ is the norm of the ideal $(\alpha)$ \cite{KONO}.  We can use this fact write our functions $F_{r, \d}(z)$ as $q$-series or linear combinations of $q$-series with complex units $\phi((\alpha))$ as coefficients over ideals $(\alpha)$ in $\Z[i]=\mathcal{O}_{\Q(i)}$ (resp. $\Z\big[\frac{1+\sqrt{-3}}{2}\big]=\mathcal{O}_{\Q(\sqrt{-3})}$).  This gives a representation of $A_r(p)$ in terms of the factorization of primes $p\equiv1\pmod{4}$ as $p=(a+bi)(a-bi)=a^2+b^2$ and primes $p\equiv1\pmod{3}$ as $p=(x+y \sqrt{-3})(x-y \sqrt{-3})=x^2+3 y^2$.  The congruences for primes $p\equiv3\pmod{4}$ (resp. primes $p\equiv2\pmod{3}$) come from Serre's result that $A_{r}(np) \equiv 0 \pmod p$ for all $n \in \Z$.  Showing that this does not happen for any prime $p \equiv 1 \pmod 4$ (resp. $p \equiv 1 \pmod 3$) will show that all of the Ramanujan-type congruences have been found for every prime except $p=2$.

Because of the way prime ideals factor in the principal ideal domains above, we will make use of the following theorem by Fermat.

\begin{thm}[Fermat]
The following are true.
\begin{enumerate}
\item Every prime $p\equiv1\pmod{4}$ has exactly one representation of the form $a^2+b^2=p$, where $a,b>0$. 
\item Every prime $p \equiv 1 \pmod 3$ has exactly one representation of the form $x^2 + 3y^2 = p$, where $x,y>0$.
\end{enumerate}
\end{thm}

\section{Proof of Theorems 1, 2, and 3}

\subsection{Combinatorial Analysis of Modular Forms with CM}

Let us define $\sum_{n=0}^\infty a_r(n)q^n = \prod_{n=1}^\infty (1-q^n)^r$.  Then, using the definition of $t$ as $t=\ell s - r$, \[\sum_{n=0}^\infty p_{-t}(n)q^n = \prod_{n=1}^\infty \frac{1}{(1-q^n)^t} = \prod_{n=1}^\infty \frac{(1-q^n)^r}{(1-q^n)^{\ell s}} = \left(\prod_{n=1}^\infty \frac{1}{(1-q^n)^{\ell s}} \right)\left(\sum_{n=0}^\infty a_r(n)q^n\right).\]
So then we have \[\left(\sum_{n=0}^\infty a_r(n)q^n\right) = \left(\prod_{n=1}^\infty {(1-q^n)^{\ell s}} \right) \left(\sum_{n=0}^\infty p_{-t}(n)q^n \right) \equiv \left(\prod_{n=1}^\infty {(1-q^{\ell n})^{s}} \right)\left(\sum_{n=0}^\infty p_{-t}(n)q^n\right) \pmod {\ell}.\] It can be seen that if $a_r(n) \equiv 0 \pmod {\ell}$ for some $n$, then $p_{-t}(n) \equiv 0 \pmod \ell$ since \\ $\prod_{n=1}^\infty (1-q^{\ell n})^{s}$ is supported only on multiples of $\ell$.  To complete the proofs of Theorems \ref{main}, \ref{two}, and \ref{three}, we will show that $a_{r} \left( \ell n + \frac{r(\ell^{2}-1)}{24} - \ell\Big\lfloor\frac{r(\ell^2-1)}{24\ell}\Big\rfloor \right) \equiv 0 \pmod \ell$ when $\ell$ is in the correct residue class.  Recall that  $F_{r, \d}(z) = \sum_{n=0}^\infty A_r(n)q^n$ and from Theorem $4$ that $F_{r, \d}(z) \vert T_{\ell} \equiv 0 \pmod {\ell}$.  This implies that
\begin{equation*}
\sum\limits_{n=0}^\infty\left(A_r(\ell n)+\chi(\ell)\ell^{\frac{r}{2}-1}A_r\left(\frac{n}{\ell}\right)\right)q^n\equiv0\pmod{\ell}.
\end{equation*}
Thus we have
\begin{equation*}
\sum\limits_{n=0}^\infty A_r(\ell n)q^n\equiv-\chi(\ell)\ell^{\frac{r}{2}-1}\sum\limits_{n=0}^\infty A_r\left(\frac{n}{\ell}\right)q^n\pmod{\ell}.
\end{equation*}
Setting $n=\d\left(n-d+\frac{r\ell}{24}\right)$, where $d:= \Big\lfloor\frac{r(\ell^2-1)}{24\ell}\Big\rfloor$, and using the relation
\begin{equation*}
\sum\limits_{n\geq0}A_r(n)q^n=\sum\limits_{n\geq0}a_r(n)q^{\d n+\frac{\d r}{24}},
\end{equation*}
we obtain
\begin{equation*}
\sum\limits_{n=0}^\infty a_r\left(\ell(n-d)+\frac{r(\ell^2-1)}{24}\right)q^n\equiv-\chi(\ell)\ell^{\frac{r}{2}-1}\sum\limits_{n=0}^\infty a_r\left(\frac{n-d}{\ell}\right)q^n\pmod{\ell}.
\end{equation*}
Replacing $n$ by $\ell k+d$ on the right hand side of the above relation, we are led to
\begin{equation*}
\sum\limits_{n=0}^\infty a_r\left(\ell(n-d)+\frac{r(\ell^2-1)}{24}\right)q^n\equiv-\chi(\ell)\ell^{\frac{r}{2}-1}q^d\prod\limits_{n=1}^\infty(1-q^{\ell n})^r\pmod{\ell}.
\end{equation*}
This proof gives a stronger result about  $a_{r} \left( \ell n + \frac{r(\ell^{2}-1)}{24} - \ell\Big\lfloor\frac{r(\ell^2-1)}{24\ell}\Big\rfloor \right)$ than is actually necessary.  

\subsection{Examples}

Recall that Olsson proved that congruence relations of this type on progressions $\ell n + a$ only occur when $a$ is in exactly the residue classes given by Theorems \ref{main}, \ref{two}, and \ref{three}.  We will now show for exactly which primes these congruence relations exist. It is noted above that $\eta^r(\d z)$ is either a newform or a linear combination  of newforms with CM.  The representation using the sum over ideals can then be used to gain insight about the coefficients of the eta-products.  The cases when $r=8$ and $r=26$ are worked out as examples below.
For $r=8$, we have \[\sum_{n=0}^\infty A_8(n) q^n = \eta^{8}(3z) = q\prod_{n=1}^\infty (1-q^{3n})^8 .\]
Using the ideal representation,
\[ \sum_{n=0}^\infty A_8(n) q^n = \sum_{(\alpha)} \phi((\alpha)) q^{N((\alpha))} = \sum_{(\alpha)} (\alpha)^3 q^{(\alpha)(\overline{\alpha})}. \]
The field is $\Z\big[\frac{1+ \sqrt{-3}}{2} \big]$ in this case, so prime ideals factor as $(p) = \left(x + y\sqrt{-3}\right)\left(x - y\sqrt{-3}\right)$.  Then we have 
\[\sum_{p} A_8(p)q^p = \sum_{(p)} (p)^3 q^{(p) (\overline{p})} = \sum_{(p)} (p)^3 q^{x^2 + 3y^2}\]
where $p$ denotes the sum over primes. We have already proven congruences for primes $p \equiv 2 \pmod 3$ for $r=8$, so now we want to study primes $p \equiv 1 \pmod 3$.  From the theorem by Fermat, we know that primes $p \equiv 1 \pmod 3$ have a unique representation as $p=x^2 + 3y^2$, so selecting the $x$ and $y$ values that produce the prime we want should also give us the coefficient of that term.  We see that
\[ \sum_{(p)} (p)^3 q^{x^2 + 3y^2} = \left( \left(x + y\sqrt{-3}\right)^3 + \left(x - y\sqrt{-3}\right)^3 \right)q^{x^2 + 3y^2},\] so the coefficients are of the form 
\[A_8(p) =  \left(x + y\sqrt{-3}\right)^3 + \left(x - y\sqrt{-3}\right)^3 = 2x(x^2 -9y^2),\]
where $p = x^2 + 3y^2$.  Requiring that $x \equiv 1 \pmod 3$ and $y > 0$ gives the correct coefficients as shown in the following table.

\begin{center}
\begin{tabular}{ | l || c | c | c | }
\hline
$p \equiv 1 \pmod 3$ & $(x, y)$ such that $p = x^2 + 3y^2$ & $A_8(p) = 2x(x^2 - 9y^2)$ \\ \hline
$7$ & $(-2, 1)$ & $20$  \\ \hline
$13$ & $(1, 2)$ & $-70$ \\ \hline
 $19$ & $(4, 1)$ & $56$ \\ \hline
$31$ & $(-2, 3)$ & $308$ \\ \hline
$37$ & $(-5, 2)$ & $110$ \\ 
\hline
\end{tabular}
\end{center}

It can be seen that none of the coefficients in the table are divisible by the corresponding primes.  In general, this can also be seen by reducing the general form of the coefficients $\mod p$ as follows.
\[A_8(p) = 2x(x^2 -9y^2) = 2x(x^2 + 3y^2) -24xy^2 \equiv -24xy^2 \pmod {p=x^2 + 3y^2}.\]
Since $p = x^2 + 3y^2$, we have that $p \nmid x$ and $p \nmid y$.  Then, since $24 = 2^3\cdot3$ and $p$ does not divide 2 or 3, we have that $-24xy^2 \not\equiv 0 \pmod p$.  This shows that there are no primes $p \equiv 1 \pmod 3$ such that $A_8(np) \equiv 0 \pmod p$ for any $n \in \Z$.  Therefore there are no simple congruence relations for those primes for $r=8$.  

The same method can be used to show that \[A_4(p) = 2x \not\equiv 0 \pmod {p = x^2 + 3y^2}\] where $\Z\big[\frac{1+ \sqrt{-3}}{2}\big]$ is the ring, $p \equiv 1 \pmod 3$, $x \equiv 1 \pmod 3$, and $y>0$, 
and that \[A_6(p) = 2(a^2 - b^2) \equiv -4b^2 \not\equiv 0 \pmod {p = a^2 + b^2}\] where $\Z[i]$ is the ring,  $p \equiv 1 \pmod 4$, $a \equiv 1 \pmod 2$, and $b >0$.

For $r = 10, 14$ and 26, the eta functions are no longer newforms but linear combinations of newforms.  The same ideas can be used, but the process becomes more complicated.  The case when $r=26$ is shown as an example below.

For $r=26$, we have \[\sum_{n=0}^\infty A_{26}(n)q^n = \eta^{26}(12z) = q^{13} \prod_{n=1}^\infty (1-q^{12n})^{26}.\]
This is a linear combination of newforms, so the ideal representation is 
\begin{eqnarray*}
&&\sum_{n=0}^\infty A_{26}(n)q^n\\
&=&\frac{1}{32617728} \Big[\sum_{(\alpha)} \phi_{+}'((\alpha))q^{N((\alpha))} + \sum_{(\alpha)} \phi_{-}'((\alpha))q^{N((\alpha))}\\
&&\hspace{2cm}- \sum_{(\alpha)} \phi_{+}''((\alpha))q^{N((\alpha))} - \sum_{(\alpha)} \phi_{-}''((\alpha))q^{N((\alpha))}\Big].
\end{eqnarray*}
The case when $r=26$ is unique (even when compared to the cases when $r=10$ and $r=14$) because it makes use of two quadratic fields; $K' = \Q(\sqrt{-3})$ and $K''=\Q(i)$.  Prime ideals that factor, factor as $(p) = \left(x + y\sqrt{-3}\right)\left(x - y\sqrt{-3}\right)$ in $\Z\big[\frac{1+\sqrt{-3}}{2}\big]$ and as $(p) = (a+bi)(a-bi)$ in $\Z[i]$.  We define $\phi_{\pm}'((\alpha)) := (-1)^{\frac{x \mp y - 1}{2}}(\alpha)^{12}$ and $\phi_{\pm}''((\alpha)) := (-1)^{3a} (\pm i)^{3b} (\alpha)^{12}$.  Using these definitions, we can write
\begin{eqnarray*}
&&\sum_{p} A_{26}(p)q^p\\
&=&\frac{1}{32617728} \Big[\sum_{(p)} \phi_{+}'((p))q^{(p)(\overline{p})} + \sum_{(p)} \phi_{-}'((p))q^{(p)(\overline{p})} - \sum_{(p)} \phi_{+}''((p))q^{(p)(\overline{p})} - \sum_{(p)} \phi_{-}''((p))q^{(p)(\overline{p})}\Big].
\end{eqnarray*}
We have already shown congruences for primes $p \equiv 11 \pmod {12}$ (which corresponds to the primes that are both $2 \pmod 3$ and $3 \pmod 4$), so now we want to study primes $p \equiv 1, 5$, or $7 \pmod {12}$ (which correspond to the primes that are $1 \pmod 3$ or $1 \pmod 4$).  From the theorem above, we know that primes $p \equiv 1 \pmod 3$ have a unique representation of the form $x^2 + 3y^2$ and primes $p \equiv 1 \pmod 4$ have a unique representation of the form $a^2 + b^2$, so selecting $x,y,a,$ and $b$ such that $x^2 + 3y^2 = a^2 + b^2 = p$ should give us the coefficients of the primes $p\equiv1 \pmod {12}$.  Primes which are $5 \pmod {12}$ (resp. $7 \pmod {12}$) have a unique representation of the form $a^2 + b^2$ (resp. $x^2 + 3y^2$) but cannot be represented by $x^2 + 3y^2$ (resp. $a^2 + b^2$), so $x$ and $y$ (resp. $a$ and $b$) are both $0$ in those cases.  Then we have the coefficients
\begin{eqnarray*}
A_{26}(p)&=&\frac{1}{32617728} \Big[ (-1)^{\frac{x-y-1}{2}} (x+y\sqrt{-3})^{12} + (-1)^{\frac{x+y-1}{2}} (x-y\sqrt{-3})^{12}\\\
 &&\hspace{2cm}+ (-1)^{\frac{x+y-1}{2}} (x+y\sqrt{-3})^{12} + (-1)^{\frac{x-y-1}{2}} (x-y\sqrt{-3})^{12}\\
&&\hspace{2cm}- (-1)^{3a}(i)^{3b} (a+bi)^{12} - (-1)^{3a} (i)^{-3b} (a-bi)^{12}\\
&&\hspace{2cm}- (-1)^{3a} (-i)^{3b} (a+bi)^{12} - (-1)^{3a} (-i)^{-3b} (a-bi)^{12} \Big].
\end{eqnarray*}
So
\begin{eqnarray*}
A_{26}(p)&=&\frac{1}{32617728} \Big[ \left( (-1)^{\frac{x-y-1}{2}} + (-1)^{\frac{x+y-1}{2}} \right)  \left(  (x+y\sqrt{-3})^{12} + (x-y\sqrt{-3})^{12} \right)\\
&&\hspace{2cm}- \left((-1)^{3a}(i)^{3b} + (-1)^{3a}(-i)^{3b} \right) \left((a+bi)^{12} + (a-bi)^{12} \right) \Big]
\end{eqnarray*}
The coefficients for the first few primes are given in the following table.

\begin{center}
\begin{tabular}{ | l || c | c | c | }
\hline
$p \equiv 1, 5, 7 \pmod {12}$ & $(a, b, x, y)$ such that $p = a^2 + b^2 = x^2 + 3y^2 $ & $A_{26}(p)$  \\ \hline
$5$ & $(-2, 1, 0, 0)$ & $0$  \\ \hline
$7$ & $(0, 0, -2, 1)$ & $0$ \\ \hline
 $13$ & $(3, 2, -1, 2)$ & $1$ \\ \hline
$17$ & $(6, 1, 0, 0)$ & $0$ \\ \hline
$19$ & $(0, 0, 4, 1)$ & $0$ \\ \hline
$29$ & $(2, 5, 0, 0)$ & $0$ \\ \hline
$31$ & $(0, 0, 2, 3)$ & $0$ \\ \hline
$37$ & $(1, 6, -5, 2)$ & $299$ \\ \hline
\hline
\end{tabular}
\end{center}

It can be seen that for primes which are $ 1 \pmod {12}$, none of the coefficients are divisible by the corresponding primes.  This can also be seen by reducing the coefficients $\mod p$ as shown below.
\[A_{26}(p) \equiv -55296x^6 y^6 - 2048a^6 b^6 \not\equiv 0 \pmod {p = a^2 + b^2 = x^2 + 3y^2}. \]
Primes which are $5 \pmod {12}$ and $7 \pmod {12}$ are visibly zero, because $A_{26}(n)$ is only supported on the progression $1 \pmod {12}$.  But in order for there to be a congruence, it must be true that $\eta^r(\d z)\vert T_{p} \equiv 0 \pmod {p}$.  This only occurs for $r=26$ when $A_{26}(np) \equiv 0 \pmod p$ for all $n$.  It is easy to show that if $A_{26}(p) = 0$, then $A_{26}(p^2) \not\equiv 0 \pmod p$, so there are no congruence relations for $r=26$ other than for primes $p \equiv 11 \pmod {12}$.  Note that there is actually one prime, namely $13$, which is not $11 \pmod {12}$ but for which congruence relations exist.  This is true because $13s - 26 = -13(2-s)$, so then we have
\begin{equation*}
\sum\limits_{n=0}^\infty p_{13(2-s)}(n)q^n = \prod\limits_{n=1}^\infty (1-q^n)^{13(2-s)} \equiv \prod\limits_{n=1}^\infty (1-q^{13n})^{2-s} \pmod {13}.
\end{equation*}
All of the coefficients, except $p_{13(2-s)}(13n)$, are visibly $0 \pmod {13}$.  The same method can be used to show that
\begin{eqnarray*}
A_{10}(p)&=&\frac{(-1)^a \left((i)^b - (-i)^b \right)}{96} \left((a+bi)^4 - (a-bi)^4 \right)\\
&\equiv&\frac{i (-1)^{a} \left((i)^b - (-i)^b \right)}{12} (a^3 b - a b^3)\\
&\not\equiv&0 \pmod p.
\end{eqnarray*}
This holds unless $p \equiv 1 \pmod {12}$, since $A_{10}(n)$ is only supported on the progression $5 \pmod {12}$.  In the cases where $A_{10}(p) = 0$, it is easy to show $A_{10}(5p) \not\equiv 0 \pmod p$.  Thus there are no congruences for $r=10$ except for primes $p \equiv 3 \pmod 4$ and $p=5$ (for the same reason that there are congruences for $p=13$ when $r=26$).  Similarly, for $r=14$, we have
 \begin{eqnarray*}
 A_{14}(p)&= &\frac{(-1)^{\frac{x-y-1}{2}} - (-1)^{\frac{x+y-1}{2}}}{720\sqrt{-3}} \left((x+y\sqrt{-3})^6 - (x-y\sqrt{-3})^6 \right)\\
&\equiv&\frac{(-1)^{\frac{x-y-1}{2}} - (-1)^{\frac{x+y-1}{2}}}{15} (4x^3 y^3)\\
&\not\equiv&0\pmod p.
\end{eqnarray*}
This holds unless $p \equiv 1 \pmod {12}$, since $A_{14}(n)$ is only supported on the progression $7 \pmod {12}$.  In the cases where $A_{14}(p) = 0$, it is easy to show $A_{14}(7p) \not\equiv 0 \pmod p$.  Thus there are no congruences for $r=14$ except for primes $p \equiv 2 \pmod 3$ and $p=7$.

\end{document}